\documentstyle[12pt]{article}

\author{Lyakhovsky V.D.\thanks{Supported in part by the
Russian Foundation for Fundamental Research, grant N 96-01-00569a}
\thanks{E-mail: lyakh!lyakh@niif.spb.su} \,
Mirolubov A.M.\\
Theoretical Department\\
Institute of Physics\\
SPb State University\\
Russia\\[7mm]
{\Large SPbU-IP-96-19}
}
\title{On the existence of deformed Lie-Poisson structures
for quantized groups
}

\begin{document}

\maketitle
\begin{abstract}
The geometrical description of deformation quantization based on quantum
duality principle makes it possible to introduce deformed Lie-Poisson
structure. It serves as a natural analogue of classical Lie bialgebra for
the case when the initial object is a quantized group. The explicit
realization of the deformed Lie-Poisson structure is a difficult problem. We
study the special class of such constructions characterized by quite a
simple form of tanjent vector fields. It is proved that in such a case it is
sufficient to find four Lie compositions that form two deformations of the
first order and four Lie bialgebras. This garantees the existence of two
families of deformed Lie-Poisson structures due to the intrinsic symmetry of
the initial compositions. The explicit example is presented.
\end{abstract}

\section{Introduction}

Quantum duality principle \cite{DRIN1},\cite{SEM1} asserts that quantization
of a Lie bialgebra $(A,A^{*})$ gives rise to a dual pair of Hopf algebras\ $
(U_p(A),U_p(A^{*}))$ or in dual terms -- $($Fun$_p(G),$Fun$_p(G^{*}))$.Here $
G$ and $G^{*}$ are the universal covering groups for $A$ and $A^{*}$
respectively. Every quantum algebra of this type can be interpreted as a
quantum group (in quantum formal series Hopf terms),
$$
U_p(A)\approx (\mbox{Fun}_p(G^{*})),
$$
and vice versa. This leads to the natural dualization of classical limits
for a given quantized object.

Thus the canonical form of deformation quantization of Lie bialgebra $
(A,A^{*})$ must be a 2-parametric family of Hopf algebras with two dual
classical limits. Within certain assumptions this family forms an analytic
variety ${\cal D}$ and the classical limits -- its boundary. The existence
of a variety ${\cal D}$ with such properties is equivalent to attributing
its members the quantum duality. The Lie bialgebra appears here in the form
of two vector fields tangent to ${\cal D}$.

>From the point of view of Lie-Poisson structures and their symmetries it can
be shown natural to consider the varieties of this type and their boundaries
entirely placed in the domain of noncommutative and noncocommutative Hopf
algebras. In this case the main property of quantum duality is preserved and
the obtained deformed vector fields can be treated as a deformed Lie-Poisson
construction with respect to a quantized group.

The paper is organized as follows. In section 2 we describe the main
features of canonical dual limits and varieties of the type ${\cal D}$. The
deformed (''lifted'') ${\cal D}_z$ versions of ${\cal D}$ are defined in
subsection 2.2 as analogs with noncommutative and noncocommutative
boundaries. In section 3 the general properties of multiplication and
comultiplication deforming functions in ${\cal D}_z$ are discussed. The
sufficient conditions for the existence of the deformation ${\cal
D\rightarrow D}_z$ are established in section 4. All the constructions
involved are explicitly illustrated by the example considered in Section 5.

\section{Deformed Poisson structures}

\subsection{Dual limits. Analytic variety ${\cal D}$}

Consider the variety ${\cal H}$ of Hopf algebras with fixed number of
generators. Its points $H\in {\cal H}$ are parameterized by the
corresponding structure constants. We must find a (smooth) curve in ${\cal H}
$ containing $U_p(A)$ and intersecting with the orbit Orb$($Fun$(G^{*}))$.
In the limit to be obtained the multiplication $m$ in $U_p(A)$ must become
Abelian. For the universal enveloping algebra $U(A)$ such a procedure is
trivially described by a contraction. The corresponding transformation of
basis $\{a_i\}$,
\begin{equation}
\label{e1}B(t):a_i\rightarrow a_i/t,
\end{equation}
leads to new structure constants $C_{ij}^{^{\prime }k}$
\begin{equation}
\label{e2}C_{ij}^{^{\prime }k}(t)=tC_{ij}^k.
\end{equation}
Algebras $U(A_t)$ form a line in Orb$U(A)$ with the limit point $U(A_0\equiv
$ Abelian$)$. Applying operators $B(t)$ to $U_p(A)$ we obtain a smooth
one-parametric curve
\begin{equation}
\label{e3}B(t)U_p(A)\equiv U_p(A_t)
\end{equation}
in the orbit Orb$(U_p(A))$ and in ${\cal H}$ -- the 2-parametric subset
\begin{equation}
\label{e4}\{U_p(A_t)\}_{p>0,t>0}\equiv {\cal D}(A,A^{*})
\end{equation}
formed by the dense set of smooth curves.

${\cal D}(A,A^{*})$ is an analytic variety with the coordinates $p$ and $t$
but the coproduct structure constants may tend to infinity when $
t\rightarrow 0$. According to the results of \cite{ETING} and \cite{RESH}
there always exists such a reparameterization of ${\cal D}$ that in new
coordinates both limits exist. In \cite{LYAKH2} it was demonstrated that for
a certain class ${\cal F}$ of quantizations the required reparameterization
is very simple.

The class ${\cal F}$ is fixed by the following conditions (see \cite{LYAKH1}
):

\begin{description}
\item[(a)]  equations
$$
(m(\mbox{id}\otimes S)\Delta )_{\downarrow V(A)}=(m(S\otimes \mbox{id}
)\Delta )_{\downarrow V(A)}=(\eta \epsilon )_{\downarrow V(A)}
$$
define uniquely the antipode $S$ of $U_p(A)$.

\item[(b)]
$$
((S_{\uparrow }\otimes \mbox{id})\Delta )_{\downarrow V(A)}=((S\otimes \mbox{
id})\Delta )_{\downarrow V(A)},
$$
$$
((\mbox{id}\otimes S_{\uparrow })\Delta )_{\downarrow V(A)}=((\mbox{id}
\otimes S)\Delta )_{\downarrow V(A)}
$$
(here $S_{\uparrow }$is a linear operator : $H\rightarrow H$ that coincides
with the antipode $S$ on $V(A)$ and is homomorphically extended to $H$ ; $
S_{\uparrow }({\bf 1})={\bf 1\ }$).
\end{description}

If the Hopf algebra $U_p(A)$ belongs to ${\cal F}$ the second classical
limit can be visualized by a simple change the parameters:
\begin{equation}
\label{e7}(p,t)\Rightarrow (h,t),\rule{1cm}{0cm}h=p/t.
\end{equation}
Then the coproduct in $U_h(A_t)$ becomes well defined in the limit $
t\rightarrow 0$ (with $h$ fixed) while the multiplication structure
constants in these new coordinates preserve their finite limit values. The
limit points correspond to Hopf algebras with Abelian multiplication. They
form an analytic curve ${\cal P}_h$. Each point of this curve can be
interpreted as an algebra of exponential coordinate functions on the group $
G_h^{*}$:\quad
\begin{equation}
\label{e13}lim_{t\rightarrow 0}\,\,U_h(A_t)\equiv U_h(A_0)\approx \mbox{Fun}
(G_h^{*})
\end{equation}
The limit $h\rightarrow 0$ describes the trivial contraction of $G_h^{*}$
into the Abelian group ${\cal AB}$.

Thus every deformation quantization of the type ${\cal F}$ can be written in
the form $U_h(A_t)$ (respectively Fun$_t(G_h^{*})$) that reveals two
canonical dual classical limits:
\begin{equation}
\label{e14}
\begin{array}{cc}
U_h(A_t) &  \\
\,_{h\rightarrow 0}\swarrow \quad \Updownarrow \quad \searrow
_{\;t\rightarrow 0} &  \\
\quad \;\;U(A_t)\quad \mbox{Fun}_h(G_t)\quad \mbox{Fun}(G_h^{*}) & \approx
U_h(
\mbox{Ab}) \\ \;\Updownarrow \quad \swarrow \qquad \qquad \searrow \quad
\Updownarrow &  \\
\quad \mbox{Fun}(G_t)\qquad \qquad \qquad \mbox{Fun}_h({\cal AB}) & \approx
U(A_h^{*})
\end{array}
\end{equation}

All the reasoning is invariant with respect to interchange $
A\rightleftharpoons A^{*}$.

\subsection{Deformed variety ${\cal D_z}$}

It was also demonstrated that for a given analytic variety ${\cal D}_{h,t}$
one can easily (and, up to equivalence, uniquely) reconstruct the
corresponding Lie bialgebra $(A,A^{*})$. Thus the Lie bialgebra and the
corresponding Lie-Poisson structure can be identified with the analytic
subvariety ${\cal D}_{h,t}$ of ${\cal H}$ whose boundaries are formed by the
contraction curves of $U(A_t)$ and Fun$(G_h^{*})$ intersecting in the common
limit -- $U($Ab$)\approx $Fun$({\cal AB})$.

To analyze the corresponding construction for quantum algebras and quantum
groups one must study the possibility to ''lift'' the whole picture (that is
the variety ${\cal D}_{h,t}$ with its boundaries) in the domain of
noncommutative and noncocommutative Hopf algebras in ${\cal H}$. Suppose
such lift exists, in this case the corresponding quaziclassical limits form
the contraction curves of quantized algebra and quantized group. They still
have a common limit, but now one finds in the intersection point nontrivial
Hopf algebra instead of Abelian and co-Abelian one.

According to the quantization ideology we must be able to treat the obtained
picture as a deformation of the initial ${\cal D}_{h,t}$. One is to
attribute an extra deformation parameter $z$ to the obtained varieties $
{\cal D}_{h,t,z}$, ${\cal D}_{h,t,0}{\cal \approx D}_{h,t}$.Thus to describe
the deformed Lie-Poisson construction for a given Lie bialgebra $(A,A^{*})$
one must find in ${\cal H}$ a tree-dimensional analytic subvariety ${\cal D}
_{h,t,z}$. From beneath it is delimited by the initial ${\cal D}_{h,t,0}$
.The other two facets -- ${\cal D}_{0,t,z}$ and ${\cal D}_{h,0,z}$ -- can
not be interpreted as canonical deformation quantizations. Parameter $z$
describes here the simultaneous change of both the multiplication and
comultiplication in ${\cal D}_{0,t,0}$and ${\cal D}_{h,0,0}$ respectively --
the situation similar to that appearing in quantum analogs of cotangent
bundle \cite{ALFAD}. Among the intersections of facets the lower two -- $
{\cal D}_{0,t,0}$ and ${\cal D}_{h,0,0}$ -- are the initial dual classical
limits while the third one -- ${\cal D}_{0,0,z}$ -- describes the trivial
contraction of the Hopf algebra that plays the role of commutative and
cocommutative one in the deformed case.

The vector fields $V_{o,t,z}$ and $W_{h,0,z}$ tangent to ${\cal D}_{h,t,z}$
and normal to ${\cal D}_{0,t,z}$ and ${\cal D}_{h,0,z}$ respectively play
here the role of deformed Poisson compositions for Fun$_z(G)$ and Fun$
_z(G^{*})$.

\section{Analytic properties of ${\cal D}_{h,t,z}$}

Let $m(t,h,z)$ and $\Delta (t,h,z)$ be the projections of multiplication and
comultiplication in $H_{h,t,z}\in {\cal D}_{h,t,z}$ on the space $A\wedge A$
and $A$ respectively.Then their Taylor expansions are
\begin{equation}
\label{e15}m(t,h,z)=\sum_{i,j,k=0;(i,j,k)\neq (0,1,0)}t^ih^jz^km_{ijk},
\end{equation}
\begin{equation}
\label{e16}\Delta (t,h,z)=\sum_{i,j,k=0;(i,j,k)\neq (1,0,0)}t^ih^jz^k\Delta
_{ijk}.
\end{equation}
Missing terms in the expansion correspond to the property of the deformation
quantization ${\cal D}_{h,t,0}{\cal \ }$where the multiplication in $U_h(A)$
and the comultiplication in Fun$_t(G)$ are untouched in the first order.
Here $m_{000}$ and $\Delta _{000}$ refer to the commutative multiplication
and primitive comultiplication respectively; in combersome expressions we
shall denote them simply by $m_0$ and $\Delta _0$.

Let $\mu $ and $\delta $ be the antisymmetrized compositions for $m$ and $
\Delta $,
\begin{equation}
\label{e17}\mu (t,h,z)=\sum_{
\begin{array}{c}
i,j,k=0 \\
(i,j,k)\neq (0,0,0),(0,j,0)
\end{array}
}t^ih^jz^k\mu _{ijk},
\end{equation}
\begin{equation}
\label{e18}\delta (t,h,z)=\sum_{
\begin{array}{c}
i,j,k=0 \\
(i,j,k)\neq (0,0,0),(i,0,0)
\end{array}
}t^ih^jz^k\delta _{ijk},
\end{equation}
The compositions $\mu _{100}$ and $\delta _{010}$ are just the $A$ and $
A_{}^{*}$ Lie multiplication and comultiplication.

Consider now the neighborhood of $H_{0,0,0}\in {\cal D}_{h,t,z}$ and write
the bialgebra properties of $H_{t,h,z}$ in terms of $\mu $ and $\delta $:
\begin{equation}
\label{e19}\delta \circ \mu =(m\otimes ^{\wedge }m)\circ (\mbox{id}\otimes
\tau \otimes \mbox{id})\circ (\Delta \otimes ^{\wedge }\Delta ),
\end{equation}
here
\begin{equation}
\label{e20}m\otimes ^{\wedge }m\equiv (m\otimes m)\circ (\mbox{id}-\tau
\otimes \tau ),
\end{equation}
\begin{equation}
\label{e21}\Delta \otimes ^{\wedge }\Delta \equiv (\mbox{id}-\tau \otimes
\tau )\circ (\Delta \otimes \Delta ).
\end{equation}
Inserting the expansions (\ref{e17},\ref{e18}) in the first nontrivial order
one finds:
\begin{equation}
\label{e22}
\begin{array}{c}
z^2\delta _{001}\circ \mu _{001}+th\delta _{010}\circ \mu _{100}+tz\delta
_{001}\circ \mu _{100}+hz\delta _{010}\circ \mu _{001}= \\
z^2(m_0\otimes \mu _{001}+\mu _{001}\otimes m_0)\circ (
\mbox{id}\otimes \tau \otimes \mbox{id})\circ (\Delta _0\otimes \delta
_{001}+\delta _{001}\otimes \Delta _0)+ \\ th(m_0\otimes \mu _{100}+\mu
_{100}\otimes m_0)\circ (
\mbox{id}\otimes \tau \otimes \mbox{id})\circ (\Delta _0\otimes \delta
_{010}+\delta _{010}\otimes \Delta _0)+ \\ tz(m_0\otimes \mu _{100}+\mu
_{100}\otimes m_0)\circ (
\mbox{id}\otimes \tau \otimes \mbox{id})\circ (\Delta _0\otimes \delta
_{001}+\delta _{001}\otimes \Delta _0)+ \\ hz(m_0\otimes \mu _{001}+\mu
_{001}\otimes m_0)\circ (\mbox{id}\otimes \tau \otimes \mbox{id})\circ
(\Delta _0\otimes \delta _{010}+\delta _{010}\otimes \Delta _0).
\end{array}
\end{equation}
Contrary to the situation with the ordinary Lie bialgebra quantization (that
is typical to ${\cal D}_{h,t,0}$) the four relations obtained from (\ref{e22}
) refer only to the facets of ${\cal D}_{h,t,z}$ . To describe the necessary
deformation into the nontrivial three-dimensional domain in ${\cal H}$ one
must consider the higher orders with much more complicated relations. Thus
one of the 3-d order equations (the $thz$ -coefficients) looks like
\begin{equation}
\label{e25}
\begin{array}{c}
\delta _{001}\circ \mu _{110}+\delta _{010}\circ \mu _{101}+\delta
_{110}\circ \mu _{001}+\delta _{011}\circ \mu _{100}= \\
(m_0\otimes \mu _{110}+\mu _{110}\otimes m_0)\circ (
\mbox{id}\otimes \tau \otimes \mbox{id})\circ (\Delta _0\otimes \delta
_{001}+\delta _{001}\otimes \Delta _0)+ \\ \left(
\begin{array}{c}
m_0\otimes \mu _{101}+\mu _{101}\otimes m_0+ \\
m_{001}^s\otimes \mu _{100}+\mu _{001}\otimes m_{100}^s+ \\
m_{100}^s\otimes \mu _{001}+\mu _{100}\otimes m_{001}^s
\end{array}
\right) \circ
\negthinspace (\mbox{id}\otimes \tau \otimes \mbox{id})\negthinspace \circ
\negthinspace (\Delta _0\otimes \delta _{010}+\delta _{010}\otimes \Delta
_0)+ \\ (m_0\otimes \mu _{100}+\mu _{100}\otimes m_0)\negthinspace \circ
\negthinspace (\mbox{id}\otimes \tau \otimes \mbox{id})\negthinspace \circ
\left(
\begin{array}{c}
\Delta _0\otimes \delta _{011}+\delta _{011}\otimes \Delta _0+ \\
\Delta _{001}^s\otimes \delta _{010}+\delta _{001}\otimes \Delta _{010}^s+
\\
\Delta _{010}^s\otimes \delta _{001}+\delta _{010}\otimes \Delta _{001}^s
\end{array}
\right) + \\
(m_0\otimes \mu _{001}+\mu _{001}\otimes m_0)\circ (\mbox{id}\otimes \tau
\otimes \mbox{id})\circ (\Delta _0\otimes \delta _{110}+\delta _{110}\otimes
\Delta _0).
\end{array}
\end{equation}
Note that symmetric parts denoted by $m_{ijk}^s$and $\Delta _{ijk}^s$ also
appear in this relation.

To be able to obtain the explicit realizations of deformed Lie-Poisson
construction one must investigate the possibility to simplify the
deformation equations imposing some restrictions on the bialgebras involved.

\section{Twice-first-order class of deformations}

Here we shall formulate the conditions that will guarantee the existence of
the nonempty ${\cal D}_{h,t,z}$ . In some sense this will also show the way
how to realize it explicitly.

{\bf Theorem}. Let

\begin{enumerate}
\item  $A$ be a Lie algebra and $A^{*}$ -- a Lie coalgebra fixed by the
compositions $\mu _{100}:V\wedge V\rightarrow V$ and $\delta
_{010}:V\rightarrow V\wedge V,$

\item  $\mu _{001}$ and $\delta _{001}$ be the deforming functions defining
first order deformations of $A$ and $A^{*}$ respectively,

\item  the following four pairs of compositions be Lie bialgebras: $(\mu
_{100},\delta _{010}),\quad $ $(\mu _{001},\delta _{010}),\quad $ $(\mu
_{100},\delta _{001}),\quad $ $(\mu _{001},\delta _{001}).$
\end{enumerate}

Then there exists in ${\cal H}$ the three-dimensional analytic subvariety $
{\cal D}_{h,t,z}\quad $ $(h,t,z\in {\tt K})$ such that

\begin{enumerate}
\item  for every fixed $\widetilde{z}\neq 0$ the corresponding
two-dimensional subvariety ${\cal D}_{h,t,\widetilde{z}}$ has the boundaries
consisting of one-dimensional ${\cal D}_{0,t,\widetilde{z}}$ and ${\cal D}
_{h,0,\widetilde{z}}$ and zero-dimensional ${\cal D}_{0,0,\widetilde{z}}$
subvarieties that can be identified with the contractions of (in general
noncommutative and noncocommutative) Hopf algebras $H_{0,t,\widetilde{z}}$
and $H_{h,0,\widetilde{z}}$ with the common limit $H_{0,0,\widetilde{z}}$
(in general noncommutative and noncocommutative),

\item  in the limit
$$
\lim _{z\rightarrow 0}{\cal D}_{h,t,z}={\cal D}_{h,t,0}
$$
the subvariety ${\cal D}_{h,t,0}$ describes the deformation quantization of
the Lie bialgebra $(\mu _{100},\delta _{010})$ with the curves ${\cal D}
_{0,t,0}$ and ${\cal D}_{h,0,0}$ presenting the canonical dual classical
limits for $H_{h,t,0}$, ${\cal D}_{0,0,0}$ coincides with the common
contraction limit of $H_{0,t,0}$ and $H_{h,t,0}\bullet $
\end{enumerate}

{\bf Proof. }According to the condition 2 for every $h,t,z,\in {\tt K}$ the
compositions
\begin{equation}
\label{e26}\widetilde{\mu }_{101}\equiv z\mu _{001}+t\mu _{100},
\end{equation}
and
\begin{equation}
\label{e27}\widetilde{\delta }_{011}\equiv z\delta _{001}+h\delta _{010}
\end{equation}
define on $V$ Lie algebra and Lie coalgebra correspondingly. It is easy to
verify that $(\widetilde{\mu }_{101},\widetilde{\delta }_{011})$ form a Lie
bialgebra. Inserting (\ref{e26},\ref{e27}) into the equation
\begin{equation}
\label{e28}
\begin{array}{c}
\widetilde{\delta }_{011}\circ \widetilde{\mu }_{101}= \\ (m_{000}\otimes
\widetilde{\mu }_{101}+\widetilde{\mu }_{101}\otimes m_{000})\circ (\mbox{id}
\otimes \tau \otimes \mbox{id})\circ (\Delta _{000}\otimes \widetilde{\delta
}_{011}+\widetilde{\delta }_{011}\otimes \Delta _{000})
\end{array}
\end{equation}
one gets the relation that splits into four equations. Each of them
describes the bialgebraic property of the corresponding pair $(\mu
_{100},\delta _{010}),\quad (\mu _{001},\delta _{010}),$ $(\mu _{100},\delta
_{001}),$ and $(\mu _{001},\delta _{001})$. The condition 3 guarantees these
properties.

Due to the result proved by P.Etingof and D.Kazhdan \cite{ETING} the Lie
bialgebra $(\widetilde{\mu }_{101},\widetilde{\delta }_{011})$ can be
quantized. The reasoning presented by N.Reshetikhin \cite{RESH} shows that
the corresponding Hopf algebras can be treated as deformation quantizations
of $(\widetilde{\mu }_{101},\widetilde{\delta }_{011})$. In \cite{LYAKH3} it
was proved that these quantized objects form in ${\cal H}$ an analytic
subvariety ${\cal D}_{x,y}\quad (x,y\in {\tt K})$ where parameters $x,y$
correspond to the trivial contractions of compositions $\widetilde{\mu }
_{101}$ and $\widetilde{\delta }_{011}$. This means that (at least when
conditions exposed in subsection 2.2 are valid) we can pass to the set of
bialgebras $\left\{ \left( x\widetilde{\mu }_{101},y\widetilde{\delta }
_{011}\right) \right\} $ that can be written as
\begin{equation}
\label{e29}x\widetilde{\mu }_{101}\Rightarrow \widetilde{\mu }_{101}^{\prime
}\equiv z^{\prime }\mu _{001}+t\mu _{100},
\end{equation}
\begin{equation}
\label{e30}y\widetilde{\delta }_{011}\Rightarrow \widetilde{\delta }
_{011}^{\prime \prime }\equiv z^{\prime \prime }\delta _{001}+h\delta _{010}.
\end{equation}
The corresponding quantizations ${\cal D}_{x,y}={\cal D}_{h,t,z^{\prime
},z^{\prime \prime }}$ form in ${\cal H}$ the four-dimen\-sional subvariety
with analytic coordinates $h,t,z^{\prime },z^{\prime \prime }\in {\tt K}$ .
Its boundaries contain the two-dimensional domains of classical objects -- $
{\cal D}_{0,t,z^{\prime },0}\ni U(\widetilde{\mu }_{101}^{\prime })$ and $
{\cal D}_{h,0,0,z^{\prime \prime }}\ni $Fun$(G(\widetilde{\delta }
_{011}^{\prime \prime }))$ , where $G(\widetilde{\delta }_{011}^{\prime
\prime })$ is the universal covering Lie group corresponding to the Lie
coalgebra $\widetilde{\delta }_{011}^{\prime \prime }$.

Now let us equalize the parameters $z^{\prime }$ and $z^{\prime \prime }$
$$
{\cal D}_{h,t,z,z}\equiv {\cal D}_{h,t,z}.
$$
Contrary to the four-dimensional picture in the three-dimensional domain $
{\cal D}_{h,t,z}$ thus obtained the two-dimensional facets ${\cal D}_{0,t,z}$
and ${\cal D}_{h,0,z}$ contain (in general) noncommutative and
noncocommutative objects. The lower facet ${\cal D}_{h,t,0}$ corresponds
obviously to the quantizations of initial Lie bialgebra $(A,A^{*})=(t\mu
_{100},h\delta _{010})$. Thus in the boundary of ${\cal D}_{h,t,z}$ only two
edge curves are classical:
$$
\begin{array}{c}
{\cal D}_{0,t,0}\ni U(t\mu _{100}), \\ {\cal D}_{h,0,0}\ni \mbox{Fun}
(G(h\delta _{010})).
\end{array}
$$
The third edge -- ${\cal D}_{0,0,z}$-- is the ''diagonal'' of the bialgebra $
(\mu _{001},\delta _{001})$ quantizations. There always exists in ${\cal D}
_{h,t,z}$ such a neighborhood of ${\cal D}_{h,t,0}$ where Hopf algebras $
H_{0,t,z}$ ,$H_{h,0,z}$ and $H_{0,0,z}$ are (in general) nonequivalent to
each other. Thus the curves ${\cal D}_{0,t,z}$ and ${\cal D}_{h,0,z}$ can be
treated as nontrivial contractions of $H_{0,t,z}$ and $H_{h,0,z}$ ( $
H_{0,0,z}$ being their common noncommutative and noncocommutative limit).$
\bullet $

We have proved that given two $\mu $- and two $\delta $- compositions (Lie
and co-Lie respectively) that form four Lie bialgebras and two first order
deformations one can always construct the parameterized (by $z$ ) set ${\cal
D}_{h,t,z}$ with the properties described in subsection 2.2.

{\bf Note. }The theorem holds also true after the interchange
$$
\begin{array}{c}
\mu _{100}\rightleftharpoons \mu _{001}, \\
\delta _{010}\rightleftharpoons \delta _{001}.
\end{array}
$$
This means that one really obtain two different three-dimensional pictures
from the four-dimensional one. The two subvarieties of ordinary deformation
quantizations ${\cal D}_{h,t,0,0}$ and ${\cal D}_{0,0,z^{\prime },z^{\prime
\prime }}$ can be treated on equal footing -- one can use the ''diagonal''
from any of them.$\bullet $

\section{Example}

The Theorem proved above can be used to construct the explicit examples of
the analytic variety ${\cal D}_{h,t,z}$ and the vector fields corresponding
to the deformed Lie-Poisson structure.

Consider the 6-dimensional vector space $V$ over ${\tt C}$ with the basis
$$
\{p_x,p_y,p_z,l_x,l_y,l_z\}.
$$
Fix the initial Lie bialgebra $(A,A^{*})$
by the compositions
\begin{equation}
\label{e31}\mu _{100}=\left\{
C_{l_yl_x}^{l_y}=1,C_{l_zl_x}^{l_z}=1,C_{p_yl_x}^{p_y}=-1,
C_{p_zl_x}^{p_z}=-1\right\} ,
\end{equation}
\begin{equation}
\label{e32}\delta _{010}=\left\{ D_{l_x}^{l_zl_y}=i\right\} .
\end{equation}
Chose the following first order deforming functions for algebras $A$ and $
A^{*}$ respectively:
\begin{equation}
\label{e33}\mu _{001}=\left\{ C_{p_zp_x}^{p_y}=i\right\} ,
\end{equation}
\begin{equation}
\label{e34}\delta _{001}=\left\{
D_{p_y}^{p_xp_y}=-1/2,D_{p_z}^{p_xp_z}=-1/2\right\} .
\end{equation}
It is easy to check that $\widetilde{\mu }_{101}^{\prime }$ and $\widetilde{
\delta }_{011}^{\prime \prime }$ (see(\ref{e29},\ref{e30})) are Lie and
co-Lie compositions for arbitrary complex parameters $\left( z^{\prime
},t\right) $ and $\left( z^{\prime \prime },h\right) $. Combined into four
pairs --$(\mu _{100},\delta _{010}),\quad (\mu _{001},\delta _{010}),\quad
(\mu _{100},\delta _{001})$ and $(\mu _{001},\delta _{001}).$ -- they form
the four Lie bialgebras mentioned in the Theorem.

The composition $\widetilde{\delta }_{011}^{\prime \prime }$ can be treated
as the direct sum of Heizenberg and $\widetilde{e}(2)$ Lie coalgebras.
This simplifies considerably the construction of quantum deformations of the
Lie bialgebra $\widetilde{(\mu }_{101}^{\prime },\widetilde{\delta }
_{011}^{\prime \prime })$. The direct sum structure and the solvability of
the coalgebra $\widetilde{\delta }_{011}^{\prime \prime }$ makes it possible
to use here the algorithm developed in \cite{LYAKH2}. The result is the
four-dimensional variety ${\cal D}_{h,t,z^{\prime },z^{\prime \prime }}$ of
Hopf algebras with the following multiplication and comultiplication :
\begin{equation}
\label{e35}
\begin{array}{c}
[p_x,p_y]=0,\qquad \qquad \qquad  \\
\,[p_z,p_x]=iz^{\prime }p_y,\qquad \qquad \; \\
\qquad [p_z,p_y]=-ih^4t^2\sinh \left( z^{\prime \prime }p_x\right) , \\
\;[l_z,l_x]=\frac t{z^{\prime \prime }h}\sinh \left( z^{\prime \prime
}hl_z\right) , \\
\qquad \quad [l_z,l_y]=i\frac th\left( \cosh \left( z^{\prime \prime
}hl_z\right) -1\right) , \\
\,[l_y,l_x]=tly,\qquad \qquad \quad \;\; \\
\qquad \qquad [p_x,l_x]=\frac t{z^{\prime \prime }}\left( \cosh \left(
z^{\prime \prime }hl_z\right) -1\right) ,\quad  \\
\,[p_y,l_x]=-\frac t2\left( 1+\cosh \left( z^{\prime \prime }hl_z\right)
\right) p_y-i
\frac{h^2t^2}{z^{\prime \prime }}\sinh \left( z^{\prime \prime }hl_z\right)
\exp \left( \frac{z^{\prime \prime }}2p_x\right) , \\ \,[p_z,l_x]=
\negthinspace - \frac t2\left( 1+\cosh \left( z^{\prime \prime }hl_z\right)
\right) p_z\negthinspace +\negthinspace h^2t\left( il_y+
\frac{z^{\prime \prime }t}2\sinh \left( z^{\prime \prime }hl_z\right)
\right) \exp \left( -\frac{z^{\prime \prime }}2p_x\right)\negthinspace , \\
\,[p_x,l_y]=it\sinh \left( z^{\prime \prime }hl_z\right) , \\
\,[p_y,l_x]=\negthinspace h^2t^2\negthinspace \left( \cosh \left( z^{\prime
\prime }hl_z\right) \exp \left(
\frac{z^{\prime \prime }}2p_x\right) \negthinspace -\negthinspace \exp
\left( -\frac{z^{\prime \prime }}2p_x\right) \right) \negthinspace
-\negthinspace i\frac{z^{\prime \prime }t}2p_y\sinh \left( z^{\prime \prime
}hl_z\right) , \\ \,[p_z,l_y]=\frac 12iz^{\prime \prime
2}h^2t^2\left( \cosh \left( z^{\prime \prime }hl_z\right) -1\right) \exp
\left( -
\frac{z^{\prime \prime }}2p_x\right) -\;\qquad \qquad \qquad \qquad
\\ \qquad \qquad \qquad \qquad \qquad -iz^{\prime \prime 2}h^2tl_x\exp
\left( -
\frac{z^{\prime \prime }}2p_x\right) -i\frac{z^{\prime \prime }t}2p_z\sinh
\left( z^{\prime \prime }hl_z\right) , \\ \,[p_x,l_z]=0,\qquad \qquad \quad
\\
\,[p_y,l_z]=0,\qquad \qquad \quad  \\
\quad \,[p_z,l_z]=2ht\sinh \left( \frac{z^{\prime \prime }}2p_x\right) ;
\end{array}
\end{equation}
\begin{equation}
\label{e36}
\begin{array}{c}
\Delta p_x=1\otimes p_x+p_x\otimes 1, \\
\Delta p_y=\exp \left( -
\frac{z^{\prime \prime }}2p_x\right) \otimes p_y+p_y\otimes \exp \left(
\frac{z^{\prime \prime }}2p_x\right) , \\ \Delta p_z=\exp \left( -
\frac{z^{\prime \prime }}2p_x\right) \otimes p_z+p_z\otimes \exp \left(
\frac{z^{\prime \prime }}2p_x\right) , \\ \Delta l_x=l_x\otimes \cosh \left(
z^{\prime \prime }hl_z\right) +1\otimes l_x-i\frac 1{z^{\prime \prime
}}l_y\otimes \sinh \left( z^{\prime \prime }hl_z\right) , \\
\Delta l_y=l_y\otimes \cosh \left( z^{\prime \prime }hl_z\right) +1\otimes
l_y+iz^{\prime \prime }l_x\otimes \sinh \left( z^{\prime \prime }hl_z\right)
, \\
\Delta l_z=l_z\otimes 1+1\otimes l_z.
\end{array}
\end{equation}
Note that the $z^{\prime }$ parameter appears only in one of the defining
relations (in the commutator $[p_z,p_x]$).

The two ''classical'' facets ${\cal D}_{0,t,z^{\prime },0}$ and ${\cal D}
_{h,0,0,z^{\prime \prime }}$ of the variety ${\cal D}_{h,t,z^{\prime
},z^{\prime \prime }}$ contain $H_{0,t,z^{\prime },0}$ -- the universal
enveloping algebra of a Lie algebra (from now on we write down only the
nonzero commutators and nontrivial cocommutators):
\begin{equation}
\label{e37}
\begin{array}{c}
[p_z,p_x]=iz^{\prime }p_y, \\
\,[l_z,l_x]=tl_z,\qquad [p_y,l_x]=-tp_y, \\
\,[l_y,l_x]=tl_y,\qquad [p_z,l_x]=-tp_z;
\end{array}
\end{equation}
and an algebra $H_{h,0,0,z^{\prime \prime }}$ of exponential coordinate
functions of the group $\widetilde{E}(2)\times E(2)$ . In the latter case
the restriction $\Delta _{\downarrow V}$ of the coproduct in $
H_{h,0,0,z^{\prime \prime }}$ coincides with that of the general case and is
presented by the relations (\ref{e36}).

The ''lower'' -- ${\cal D}_{h,t,0,0}$ --and the ''upper'' -- ${\cal D}
_{0,0,z^{\prime },z^{\prime \prime }}$ -- facets of ${\cal D}_{h,t,z^{\prime
},z^{\prime \prime }}$ look like the ordinary deformation quantizations. The
points $H_{h,t,0,0}$ refer to the quantizations of $\left( \mu _{100},\delta
_{010}\right) $-bialgebra with the defining relations
\begin{equation}
\label{e38}
\begin{array}{c}
\,[l_z,l_x]=tl_z,\qquad [p_y,l_x]=-tp_y-ih^3t^2l_z, \\
\,[l_y,l_x]=tl_y,\qquad [p_z,l_x]=-tp_z+ih^2tl_y;\;
\end{array}
\end{equation}
\begin{equation}
\label{e39}\Delta l_x=l_x\otimes 1+1\otimes l_x-ihl_y\otimes l_z.
\end{equation}
While the points $H_{0,0,z^{\prime },z^{\prime \prime }}$ describe the
quantized Heizenberg algebra
\begin{equation}
\label{e40}[p_z,p_x]=iz^{\prime }p_y,
\end{equation}
\begin{equation}
\label{e41}
\begin{array}{c}
\Delta p_y=\exp \left( -
\frac{z^{\prime \prime }}2p_x\right) \otimes p_y+p_y\otimes \exp \left(
\frac{z^{\prime \prime }}2p_x\right) , \\ \Delta p_z=\exp \left( -\frac{
z^{\prime \prime }}2p_x\right) \otimes p_z+p_z\otimes \exp \left( \frac{
z^{\prime \prime }}2p_x\right) .
\end{array}
\end{equation}

Now let us put $z^{\prime }=z^{\prime \prime }=z$ and consider the
three-dimensional variety ${\cal D}_{h,t,z}$ .This means that we chose only
the ''diagonal'' points in $H_{0,0,z^{\prime },z^{\prime \prime }}$ . The
lower facet $H_{h,t,0,0}$ rests unchanged and fixes the initial (in our
interpretation the undeformed) Lie-Poisson structure -- the Heizenberg
Poisson algebra for the group $G(\mu _{100})$ or the $\mu _{100}$ Poisson
algebra for the Heizenberg group.

The two-dimensional facets ${\cal D}_{0,t,z}$ and ${\cal D}_{h,0,z}$ can no
more be considered as ordinary quantizations. Each of them have only one
classical edge ${\cal D}_{0,t,0}\in {\cal D}_{0,t,z}$ and ${\cal D}
_{h,0,0}\in {\cal D}_{h,0,z}$, their common edge ${\cal D}_{0,0,z}$ refers
to the quantized Heizenberg algebra $H_{0,0,z}$ (defined by (\ref{e40},\ref
{e41}) with $z^{\prime }=z^{\prime \prime }=z$ ). Thus considering the
parameter $z$ in the variety ${\cal D}_{h,0,z}$ as the deformation parameter
we see that in the Hopf algebras $H_{h,0,z}$ ,
\begin{equation}
\label{e42}
\begin{array}{c}
[p_z,p_x]=izp_y, \\
\Delta p_y=\exp \left( -\frac z2p_x\right) \otimes p_y+p_y\otimes \exp
\left( \frac z2p_x\right) , \\
\Delta p_z=\exp \left( -\frac z2p_x\right) \otimes p_z+p_z\otimes \exp
\left( \frac z2p_x\right) , \\
\Delta l_x=l_x\otimes \cosh \left( zhl_z\right) +1\otimes l_x-i\frac
1zl_y\otimes \sinh \left( zhl_z\right) , \\
\Delta l_y=l_y\otimes \cosh \left( zhl_z\right) +1\otimes l_y+izl_x\otimes
\sinh \left( zhl_z\right) ,
\end{array}
\end{equation}
not only the multiplication but also the comultiplication is deformed. The
same is true also for $H_{0,t,z}$ :
\begin{equation}
\label{e43}
\begin{array}{c}
[p_z,p_x]=izp_y, \\
\,[l_z,l_x]=tl_z,\qquad [p_y,l_x]=-tp_y, \\
\,[l_y,l_x]=tl_y,\qquad [p_z,l_x]=-tp_z; \\
\Delta p_y=\exp \left( -\frac z2p_x\right) \otimes p_y+p_y\otimes \exp
\left( \frac z2p_x\right) , \\
\Delta p_z=\exp \left( -\frac z2p_x\right) \otimes p_z+p_z\otimes \exp
\left( \frac z2p_x\right) ,
\end{array}
\end{equation}
This property of ${\cal D}_{h,0,z}$ and ${\cal D}_{0,t,z}$ resembles that of
the quantized version of the cotangent bundle \cite{ALFAD} where the
multiplication on Fun$(T^{*}G)$ is deformed simultaneously with the group $
G\rhd T^{*}$ itself. Note that in our case the Hopf structure is preserved.

The varieties ${\cal D}_{h,t,\widetilde{z}}$ with fixed $z=\widetilde{z}\neq
0$ are especially interesting for us. Their boundaries lie in ${\cal D}_{h,0,
\widetilde{z}}$ and ${\cal D}_{0,t,\widetilde{z}}$ . The internal points $
H_{h,t,\widetilde{z}}$ can be interpreted as quantizations of $H_{0,t,
\widetilde{z}}$ by $H_{h,0,\widetilde{z}}$ and $\widetilde{z}$ measures the
difference between the canonical -- ${\cal D}_{h,t,0}$ -- and deformed -- $
{\cal D}_{h,t,\widetilde{z}}$ -- pictures. It is evident that Hopf algebras $
H_{0,t,\widetilde{z}}$ (see (\ref{e40},\ref{e41}) are inequivalent to the
quantized Heizenberg algebra $H_{0,0,\widetilde{z}}$ . The same is true for $
H_{h,0,\widetilde{z}}$ and $H_{0,0,\widetilde{z}}$. Thus the boundaries $
{\cal D}_{h,0,\widetilde{z}}$ and ${\cal D}_{0,t,\widetilde{z}}$ can be
treated as the nontrivial contraction curves with the common contraction
limit ${\cal D}_{0,0,\widetilde{z}}$ .

In the initial variety ${\cal D}_{h,t,0}$ the Lie-Poisson structure is
defined by the tangent fields $V_{0,t,0}$:
\begin{equation}
\label{e44}\delta (l_x)=-il_y\wedge l_z,
\end{equation}
and $W_{h,0,0}$ :
\begin{equation}
\label{e45}
\begin{array}{c}
\mu (l_z,l_x)=l_z,\quad \qquad \qquad \mu (l_y,l_x)=l_y, \\
\mu (p_y,l_x)=-p_y,\quad \mu (p_z,l_x)=-p_z+ih^2l_y.
\end{array}
\end{equation}
In the deformed case the field $V_{0,t,\widetilde{z}}$ contains both
multiplication and comultiplication components. To simplify the comparison
with $V_{0,t,0}$ we shall present the first terms of its power series
expansion in $\widetilde{z}$ and $t$ :
\begin{equation}
\label{e46}
\begin{array}{c}
\mu (l_z,l_y)=\frac i2t
\widetilde{z}^2l_z^2+..., \\ \mu (p_x,l_y)=it
\widetilde{z}l_z^{}+..., \\ \mu (p_y,l_y)=-\frac i2t
\widetilde{z}^2l_z^{}p_y+..., \\ \mu (p_z,l_y)=-\frac i2t
\widetilde{z}^2l_z^{}p_z+..., \\ \mu (p_z,l_z)=t
\widetilde{z}p_x+..., \\ \delta (l_x)=-il_y\wedge l_z+..., \\
\delta (l_y)=i\widetilde{z}^2l_x\wedge l_z+...\quad .
\end{array}
\end{equation}
The field $V_{0,t,0}$ is reobtained in the limit $\lim _{\widetilde{z}
\rightarrow 0}V_{0,t,\widetilde{z}}=V_{0,t,0}$ . The co-Poisson structure
defined by $V_{0,t,0}$ is deformed to describe the correlation between the
possible additional quantization of $H_{0,t,\widetilde{z}}$ and the
deformation of its multiplication structure constants.

The field $W_{h,0,\widetilde{z}}$ does not obtain the additional components
in the coproduct sector:
\begin{equation}
\label{e47}
\begin{array}{c}
\mu (l_z,l_x)=l_z+...,\quad \mu (l_y,l_x)=l_y,\quad \mu (l_z,l_y)=\frac i2h
\widetilde{z}^2l_z^2+..., \\ \mu (p_y,l_x)=-p_y+...,\qquad \quad \qquad \mu
(p_z,l_x)=-p_z+ih^2l_y+..., \\
\mu (p_x,l_x)=\frac 12h^2
\widetilde{z}l_z^2+...,\qquad \qquad \mu (p_x,l_y)=ih\widetilde{z}
l_z+...,\qquad \; \\ \mu (p_y,l_y)=-\frac i2h
\widetilde{z}^2l_zp_y+...,\;\;\qquad \mu (p_z,l_y)=-\frac i2h\widetilde{z}
^2l_zp_z+..., \\ \mu (p_z,l_z)=h\widetilde{z}p_x+...\quad .
\end{array}
\end{equation}
Considered as function of $z$ the field $W_{h,0,z}$ describes the explicit
dependence of the Poisson algebra on the quantization of the group $H_{0,t,z}
$ where $z$ now plays the role of the deformation parameter.

\section{Acknowledgments}

The authors express their gratitude to Prof. P.P.Kulish for fruitful
discussions.

\end{document}